\documentclass[reqno]{amsart}
\usepackage{amssymb,latexsym,amscd,graphicx}
\newtheorem{lem}{Lemma}[section]
\newtheorem{cor}{Corollary}[section]
\newtheorem{tem}{Theorem}[section]

\theoremstyle{definition}
\newtheorem{defn}{Definition}[section]
\theoremstyle{remark}
\newtheorem{rem}{Remark}[section]
\begin{document}
\title{Hypergeometric Zeta Functions} 
\author{Abdul Hassen and Hieu D. Nguyen}
\date{7/10/05}
\address{Department of Mathematics, Rowan University, Glassboro, NJ 08028.}
\email{hassen@rowan.edu, nguyen@rowan.edu}
\subjclass[2000]{Primary 11M41}
\keywords{Riemann zeta function, hypergeometric integrals}
\begin{abstract}
This paper investigates a new family of special functions referred to as hypergeometric zeta functions.  Derived from the integral representation of the classical Riemann zeta function, hypergeometric zeta functions exhibit many properties analogous to their classical counterpart, including the intimate connection to Bernoulli numbers.  These new properties are treated in detail and are used to demonstrate a functional inequality satisfied by second-order hypergeometric zeta functions.

\end{abstract}

\maketitle
\def\theequation{\thesection.\arabic{equation}}
\section{INTRODUCTION}
\setcounter{equation}{0}
Riemann demonstrated in \cite{R} that the classical zeta function

$$\zeta(s)=\sum_{n=1}^{\infty}\frac{1}{n^s}$$
admits the integral representation

\begin{equation} \zeta(s)=\frac{1}{\Gamma(s)}\int_0^\infty \frac{x^{s-1}}{e^x-1}\,dx. \label{eq 1.1}\end{equation}

By using complex analysis he was able to relate  $\left(\ref{eq 1.1}\right)$ to a suitable contour integral that allowed him to continue $\zeta(s)$ analytically to the entire complex plane (except for a simple pole at $s=1$  ) and to establish its celebrated functional equation: 
\begin{equation} \zeta(s)=2 (2 \pi)^{s-1} \sin \left( \frac{\pi}{2}s\right)\Gamma(1-s)\zeta(1-s). \label{eq 1.2}\end{equation}

Riemann's proof of $\left(\ref{eq 1.2}\right)$ (he actually gave two proofs) used residue theory, an effective strategy here since the integrand in $\left(\ref{eq 1.1}\right)$  has nice singularities on the complex plane.  In particular, the roots of $e^z-1=0$  are located at integer multiplies of $2\pi i$, i.e. $z=2\pi i n$, and allows for an easy calculation of the corresponding residues.

In this paper, we investigate an interesting generalization of $\left(\ref{eq 1.1}\right)$ that fleshes out the important role acted out by its singularities.  To this end, we replace the denominator $e^x-1$  in $\left(\ref{eq 1.1}\right)$ by  an arbitrary Taylor difference  $e^x-T_{N-1}(x)$, where $N$ is a positive integer and $T_{N-1}(x)$  is the Taylor polynomial of $e^x$  at the origin having degree $N-1$.  This defines a family of higher-order zeta functions denoted by:\begin{equation} \zeta_N(s)=\frac{1}{\Gamma(s+N-1)}\int_0^\infty \frac{x^{s+N-2}}{e^x-T_{N-1}(x)}\,dx\;\;\;\;\;\;(N\geq 1). \label{eq 1.3}\end{equation}  Observe that $\zeta_1(s)=\zeta(s)$.  For reasons to be explained later, we shall refer to $\{\zeta_N(s)\}$ as {\it hypergeometric zeta functions}.  Following Riemann, we develop their analytic continuation to the entire complex plane, except for $N$ simple poles at $s=1,0,-1, \cdots, 2-N$, and establish many properties analogous to those satisfied by Riemann's zeta function.\\

A classical property of $\zeta(s)$ is its evaluation at negative integers.  Euler demonstrated that its values are expressible in terms of Bernoulli numbers:
$$\zeta(-n)=-\frac{B_{n+1}}{n+1}.$$
Here, the Bernoulli numbers $B_n$  are generated by
$$\frac{x}{e^x-1}=\sum_{n=0}^{\infty}\frac{B_n}{n!} x^n.$$ 

In the case of hypergeometric zeta functions given by $\left(\ref{eq 1.3}\right)$, we find that they can be similarly expressed in terms of generalized Bernoulli numbers.  For example, when $N=2$  it is shown that$$\zeta_2(-n)=(-1)^{n+1}\frac{ 2 B_{2,n+1}}{n(n+1)}.$$ The coefficients $B_{2,n}$  above are likewise generated by $$\frac{x^2/2}{e^x-1-x}=\sum_{n=0}^{\infty}\frac{B_{2,n}}{n!} x^n.$$ F.T. Howard initiated a study of these coefficients in \cite{H1}-\cite{H3}.  He referred to them as $A_n$  and discovered many interesting properties analogous to those of the classical Bernoulli numbers.  In particular, Howard used Hadamard Factorization Theorem to express these numbers as \begin{equation} B_{2,n}=-n!\sum_{k=1}^\infty r_k^{-n}\cos \left( n \theta_k \right). \label{eq 1.4}\end{equation} Here, $z_k =x_k+i y_k=r_ke^{i\theta_k}$ are the zeros of $e^z-1-z=0$  that are located in the upper-half plane.  As a result, Howard in \cite  {H1} established the inequality \begin{equation} |B_{2,n}|<\frac{n!}{(2\pi)^n}\zeta(n). \label{eq 1.5}\end{equation} 

In this present work, we extend $\left(\ref{eq 1.4}\right)$  by interpreting it as a discrete case of the following continuous result, valid for $\Re(s)<0$ (see Theorem 4.1):\begin{equation} \zeta_2(s)=-2\Gamma(-s)\sum_{k=1}^\infty r_k^{s-1}\cos \left[(s-1)(\pi-\theta_k)\right]. \nonumber \end{equation} Using Howard's estimate for the size of the roots $z_k$, we obtain an inequality between $\zeta_2(s)$  and $\zeta(s)$  that generalizes $\left(\ref{eq 1.5}\right)$  for  $\Re(s)<0$ (see Theorem 4.2):  \begin{equation}\left| \zeta_2(s)\right|<2(2\pi)^{\Re(s)}\left|\Gamma(-s)\right | e^{\Im(s) (\pi-\theta_1)}\zeta(1-\Re(s)), \label{eq 1.6} \end{equation} where $\theta_1 \approx 1.2978$ is the angle of the smallest nonzero root of $e^z-1-z=0$ in the upper half-plane. Since $\zeta(1-\Re(s))< \zeta_2(1-\Re(s))$  (see $\left(\ref{eq 2.10}\right)$), this yields a 'functional inequality' for $\zeta_2(s)$: \begin{equation}\left| \zeta_2(s)\right|<2(2\pi)^{\Re(s)}\left|\Gamma(-s)\right | e^{\Im(s)(\pi-\theta_1)}\zeta_2(1-\Re(s)) . \label{eq 1.7} \end{equation} Observe that inequality  $\left(\ref{eq 1.7}\right)$  resembles the functional equation for $\zeta(s)$   given by $\left(\ref{eq 1.2}\right)$. The more difficult problem of course is to extend this functional inequality to an equality, which most likely will require knowing the precise locations of the zeros $\{z_k\}$ of $e^z-1-z=0$.  Some results describing the approximate location of these roots appear in Appendix I.\\

Our paper is organized as follows.  In section 2, we define hypergeometric geometric functions, establish convergence on a right half-plane, and develop their series representations.  In section 3, we reveal their analytic continuation to the entire complex plane, except at a finite number of poles, and calculate their residues in terms in generalized Bernoulli numbers.  In section 4, we establish a series formula valid on a left half-plane and use it to obtain a functional inequality satisfied by second-order hypergeometric zeta functions and to prove a conjecture made by Howard in \cite{H3} regarding the growth of generalized Bernoulli numbers.  Sections 5 and 6 are appendices demonstrating some results that are used in the main body of the paper regarding the zeros of $e^z-T_{N-1}(z)=0$ (Appendix I) and listing the first ten of these zeros for the cases $N=2$ and $N=3$ (Appendix II).\\

\noindent{\em Acknowledgement:} Both authors would like to thank their colleague and friend Thomas J. Osler for the many useful conversations on hypergeometric zeta functions.
 
\section{PRELIMINARIES}
\setcounter{equation}{0}
In this section we formally define hypergeometric zeta functions, establish a domain of convergence, and demonstrate their series representations.\begin{defn} Denote the Maclaurin (Taylor) polynomial of the exponential function $e^x$  by \begin{eqnarray*}T_{N}(x)=\sum_{k=0}^N \frac{x^k}{k!}.\end{eqnarray*}  We define the {\em $N^{th}$-order hypergeometric zeta function} (or just {\em hypergeometric zeta function} for short) to be
\begin{equation} \zeta_N(s)=\frac{1}{\Gamma(s+N-1)}\int_0^\infty \frac{x^{s+N-2}}{e^x-T_{N-1}(x)}\,dx\;\;\;\;\;\;(N\geq 1). \label{eq 2.1}\end{equation} Moreover, for $ N = 0$, we set $\zeta_0(s)=1$.\end{defn} \begin{rem} Observe that the first-order hypergeometric zeta function reduces to Riemann's zeta function, i.e. $\zeta_1(s)=\zeta(s)$.\end{rem}
\begin{lem}$ \zeta_N(s)$  converges absolutely for $\sigma=\Re(s)>1$.\end{lem}

\begin{proof}
Let $K>0$  be such that $e^x\geq e^{x/2}+T_{N-1}(x)$   for all $x\geq K$.  This is equivalent to $e^x-T_{N-1}(x)\geq e^{x/2}$.  For $\sigma >1$, we have
\begin{eqnarray*}\left| \zeta_N(s)\right| & \leq &   \frac{1}{|\Gamma(s+N-1)|}\left[\int_0^K \left |\frac{ x^{s+N-2} }{ e^x - T_{N-1} (x) } \right| dx + \int_K^\infty \left |\frac{x^{s+N-2}}{e^x-T_{N-1}(x)}\right|dx \right] \\
&  \leq &  \frac{1}{|\Gamma(s+N-1)|}\left[\int_0^K \frac{x^{\sigma +N-2}}{x^N/N!} dx\,+\,\int_K^\infty x^{\sigma +N-2}e^{-x/2}dx \right]\\
&  \leq &  \frac{1}{|\Gamma(s+N-1)|}\left[N!\int_0^K x^{\sigma -2} dx\,+\,2^{\sigma + N-2}\int_K^\infty y^{\sigma +N-2}e^{-y}dy \right]\\
 &  \leq &  \frac{1}{|\Gamma(s+N-1)|}\left[\frac{N!K^{\sigma-1}}{\sigma-1}\,+\,2^{\sigma + N-2}\Gamma(\sigma+N-2) \right]\\ & < & \infty.
\end{eqnarray*}
This proves our lemma.\end{proof}

The next two lemmas provide hypergeometric zeta with a series representation, which reduces formally to the harmonic series at $s=1$.
\begin{lem}For $\sigma >1$, we have
\begin{equation} \zeta_N(s)=\sum_{n=1}^\infty  f_n(N,s),\label{eq 2.2}\end{equation} where
\begin{equation} f_n(N,s)=\frac{1}{\Gamma(s+N-1)}\int_0^\infty x^{s+N-2}T_{N-1}^{n-1}(x)e^{-nx}\,dx.\;\;\;\;\;\; \label{eq 2.3}\end{equation} 
\end{lem}

\begin{proof} Since  $\left|T_{N-1}(x)e^{-x}\right|<1$ for all $x>0$, we can rewrite the integrand in $\left(\ref{eq 2.1}\right)$   as a geometric series:
$$\frac{ x^{s+N-2} }{ e^x - T_{N-1} (x) }=\frac{ e^{-x}x^{s+N-2} }{1-  T_{N-1} (x)e^{-x} }= e^{-x}x^{s+N-2}\sum_{n=0}^\infty\left[T_{N-1}(x)e^{-x}\right]^n=x^{s+N-2}\sum_{n=1}^\infty T_{N-1}^{n-1}(x)e^{-nx}.$$
The lemma now follows by reversing the order of integration and summation because of Dominated Convergence Theorem:\begin{eqnarray*}\zeta_N(s)&=&\frac{1}{\Gamma(s+N-1)}\int_0^\infty x^{s+N-2}\sum_{n=1}^\infty T_{N-1}^{n-1}(x)e^{-nx}\,dx\\&=& \sum_{n=1}^\infty\left[\frac{1}{\Gamma(s+N-1)}\int_0^\infty x^{s+N-2} T_{N-1}^{n-1}(x)e^{-nx}\,dx \right ]. \end{eqnarray*}
\end{proof}

\begin{rem}  Observe that $\zeta_2(s)$  can be expressed in terms of hypergeometric functions:
\begin{eqnarray}\zeta_2(s)&=&\frac{1}{\Gamma(s+1)}\int_0^\infty x^{s}\sum_{n=1}^\infty (1+x)^{n-1}e^{-nx}\,dx\nonumber\\
&=& \sum_{n=1}^\infty\left[\frac{1}{\Gamma(s+1)}\int_0^\infty x^{s}(1+x)^{n-1}e^{-nx}\,dx \right ]\nonumber\\
&=& \sum_{n=1}^\infty U(s+1,s+1+n,n),\label{eq 2.4} \end{eqnarray}
where $U(a,b,z) $ is the confluent hypergeometric function of the second kind defined by
\begin{eqnarray}U(a,b,z)= \frac{1}{\Gamma(a)}\int_0^\infty x^{a-1}(1+x)^{b-a-1}e^{-z x}\,dx. \label{eq 2.5}\end{eqnarray}
This justifies our use of the term `hypergeometric zeta function' for $\zeta_2(s)$.  Actually, a much more evident reason for this nomenclature in the general case can be seen directly from definition $\left(\ref{eq 2.1}\right)$, where the integrand can be expressed in terms of the confluent hypergeometric series:
$$\zeta_N(s)=\frac{\Gamma(N)}{\Gamma(s+N-1)}\int_0^\infty\frac{x^{s-1}}{_1F_1(1,N;x)-1} dx.$$
This representation is discussed further in our concluding remarks at the end of section 4.\end{rem}
\begin{lem} For $f_n(N,s)$ given by  $\left(\ref{eq 2.3}\right)$, we have 
\begin{eqnarray}f_n(N,1)=\frac{1}{n}.\label{eq 2.6}
\end{eqnarray}\end{lem}
\begin{proof} Since $x^{N-1}=(N-1)!\left[T_{N-1}(x)-T_{N-2}(x)\right]$, it follows that
\begin{eqnarray*}f_n(N,1)&=&\frac{1}{(N-1)!}\int_0^\infty x^{N-1}T_{N-1}^{n-1}(x)e^{-nx}\,dx\\
&=& \int_0^\infty T_{N-1}^{n}(x)e^{-nx}\,dx-\;\int_0^\infty T_{N-2}(x)T_{N-1}^{n-1}(x)e^{-nx}\, dx. \end{eqnarray*}
But the two integrals above merely differ by 1/n, which results from integrating by parts:
$$\int_0^\infty T_{N-1}^{n}(x)e^{-nx}\,dx=\frac{1}{n}+\int_0^\infty T_{N-2}(x)T_{N-1}^{n-1}(x)e^{-nx}\,dx.$$ This establishes the lemma.
\end{proof}

\begin{rem} We deduce from $\left(\ref{eq 2.2}\right)$ and  $\left(\ref{eq 2.6}\right)$ that $\zeta_N(1)=\sum_{n=1}^\infty 1/n$  formally generates the harmonic series for all $N$.  This reveals our motivation for normalizing the gamma factor in $\left(\ref{eq 2.1}\right)$ as we did in defining $\zeta_N(s)$.\end{rem}

To demonstrate next that $\zeta_N(\sigma)>\zeta(\sigma) $ for $\sigma >1$ and $N>1$, we shall need the help of two additional lemmas.
\begin{lem} For $ \Re(s)=\sigma>1$, we have 
\begin{equation} \zeta_N(s)=\sum_{n=1}^\infty \frac{\mu_N(n,s)}{n^{s+N-1} },\label{eq 2.7} \end{equation}
where
\begin{equation} \mu_N(n,s)=\sum_{k=0}^{(N-1)(n-1)} \frac{a_k(N,n)}{n^{k}}(s+N-1)_k.\label{eq 2.8} \end{equation}
Here $a_k(N,n)$ is generated by
$$\left(T_{N-1}(x)\right)^{n-1}=\left(\sum_{k=0}^{N-1}\frac{x^k}{k!}\right)^{n-1}=\sum_{k=0}^{(N-1)(n-1)} a_k(N,n)x^k.$$
\end{lem}
\begin{proof} With $a_k(N,n)$ as given above, we have \begin{eqnarray*}\zeta_N(s)&=&\sum_{n=1}^\infty \frac{1}{\Gamma(s+N-1)}\int_0^\infty x^{s+N-2} T_{N-1}^{n-1}(x)e^{-nx}\,dx\\
&=& \sum_{n=1}^\infty\left[\frac{1}{\Gamma(s+N-1)}\int_0^\infty x^{s+N-2}\left( \sum_{k=0}^{(N-1)(n-1)} a_k(N,n)x^k \right)e^{-nx}\,dx \right ]\\
& & \\
&=& \sum_{n=1}^\infty\left[\frac{1}{\Gamma(s+N-1)}\frac{1}{n^{s+N-1}}\int_0^\infty  \left(\sum_{k=0}^{(N-1)(n-1)} a_k(N,n)\frac{x^{s+k+N-2}}{n^k}\right)e^{-x} \,dx \right ]\\
& &\\
&=& \sum_{n=1}^\infty\left[\frac{1}{\Gamma(s+N-1)}\frac{1}{n^{s+N-1}} \left(\sum_{k=0}^{(N-1)(n-1)} \frac{a_k(N,n)}{n^k}\int_0^\infty x^{s+k+N-2} e^{-x} \,dx \right) \right ]\\
& &\\
&=&\sum_{n=1}^\infty\left[\frac{1}{n^{s+N-1}} \left(\sum_{k=0}^{(N-1)(n-1)} \frac{a_k(N,n)}{n^k}\frac{\Gamma(s+N+k-1)}{\Gamma(s+N-1)}\right)\right]\\
& &\\
&=&\sum_{n=1}^\infty\left[\frac{1}{n^{s+N-1}} \left(\sum_{k=0}^{(N-1)(n-1)} \frac{a_k(N,n)}{n^k}(s+N-1)_k\right)\right]\\
&=& \sum_{n=1}^\infty \frac{\mu_N(n,s)}{n^{s+N-1} }.\end{eqnarray*}\end{proof}
  
\begin{lem}\begin{equation} \mu_N(n,1)=\sum_{k=0}^{(N-1)(n-1)} \frac{a_k(N,n)}{n^{k}}(N)_k=n^{N-1}\label{eq 2.9} \end{equation}\end{lem}
\begin{proof}  Since $\mu_N(n,s)/n^{s+N-1} =f_n(N,s)$, we have from  $\left(\ref{eq 2.6}\right)$ that $\mu_N(n,1)/n^{N}=f_n(N,1)=1/n$. The result of the lemma now becomes clear.
\end{proof}

\begin{tem}For $N>1$ and real values of $s=\sigma >1$, we have
\begin{equation}  \zeta_N(\sigma)>\zeta(\sigma).\label{eq 2.10} \end{equation}\end{tem}
\begin{proof} It is clear from  $\left(\ref{eq 2.8}\right)$  that $\mu_N(n,s)$   is a strictly increasing function when taking on real values of $s$ since it is a polynomial with positive coefficients.  Hence, for $\sigma >1$,$$\zeta_N(\sigma)=\sum_{n=1}^\infty\frac{\mu_N(n,\sigma)}{n^{\sigma+n-1}}>\sum_{n=1}^\infty\frac{\mu_N(n,1)}{n^{\sigma+n-1}}=\sum_{n=1}^\infty\frac{1}{n^{\sigma}}=\zeta(\sigma).$$
\end{proof}

\begin{figure} \centering
\includegraphics[width=.35\textwidth]{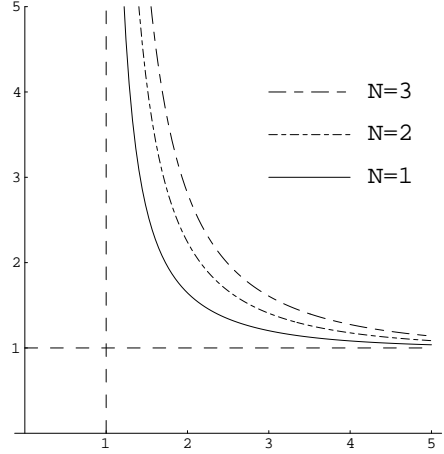}
\caption{Graphs of $\zeta_N(\sigma)$ for $N=1,2,3$.}
\label{figure:monotonicity}
\end{figure}

\begin{rem}\
\begin{enumerate}
\item[(a)] Observe that the coefficient $\mu_N(n,s)$  in the series representation of  $\zeta_N(s)$ depends on both $n$ and $s$.  In this sense it is a generalized Dirichlet series. Of course, we would like to find an expression of  $\mu_N(n,s)$ that allows us to write  $\zeta_N(s)$ as an ordinary Dirichlet series.  At the present moment, we do not know even for $N=2$ if the series representation  $\left(\ref{eq 2.4}\right)$  involving the confluent hypergeometric function  $\left(\ref{eq 2.5}\right)$  will lead to any such result.
\item[(b)] Graphical evidence (cf. Figure \ref{figure:monotonicity}) suggests the following `monotonicity' conjecture:$$\zeta_N(\sigma)>\zeta_{N-1}(\sigma).$$ \end{enumerate}
\end{rem}

\section{ ANALYTIC CONTINUATION}

\setcounter{equation}{0}      
                          
In this section we develop the analytic continuation of hypergeometric zeta to the entire complex plane.  We shall discuss two different approaches.  The first involves rewriting the integral $\left(\ref{eq 2.1}\right)$ in stages to extend the domain of  $\zeta_N(s)$ strip by strip and the second uses contour integration to perform the analytic continuation in one stroke.  As we will see each method has its advantages.\\

Assume $\Re(s)>1$.  Then $\left(\ref{eq 2.1}\right)$  can be rewritten as 
\begin{eqnarray} \Gamma(s+N-1)\zeta_N(s)&=&\int_0^1 \frac{x^{s+N-2}}{e^x-T_{N-1}(x)}\,dx\;+\;\int_1^\infty \frac{x^{s+N-2}}{e^x-T_{N-1}(x)}\,dx \label{eq 3.1}\\
  &=&\int_0^1 \left(\frac{1}{e^x-T_{N-1}(x)}-\frac{N!}{x^N}\right)x^{s+N-2}\,dx\;+\;\frac{N!}{s-1}\;+\;\int_1^\infty \frac{x^{s+N-2}}{e^x-T_{N-1}(x)}\,dx.\nonumber\end{eqnarray}
The last formula in  $\left(\ref{eq 3.1}\right)$ is analytic in the strip $0<\Re(s)\leq 1$,  except for the pole at $s=1$, since both integrals on the right hand side are convergent on this domain.  Moreover, for $0<\Re(s)<1$, $$\frac{N!}{s-1}=-\int_1^\infty\frac{x^{s+N-2}}{x^N} dx.$$ 
Therefore, we obtain the following result:
\begin{tem}For $0<\Re(s)<1$,
$$\zeta_N(s)=\frac{1}{\Gamma(s+N-1)}\int_0^\infty \left(\frac{1}{e^x-T_{N-1}(x)}-\frac{N!}{x^N}\right)x^{s+N-2}\,dx.$$
\end{tem}

\begin{rem} This process can be repeated to extend $\zeta_N(s)$ analytically  to  $-1<\Re(s)<0$, thus skipping over the second pole at $s=0$:
\begin{eqnarray*}\Gamma(s+N-1)\zeta_N(s)&=& \int_0^1 \left(\frac{1}{e^x-T_{N-1}(x)}-\frac{N!}{x^N}\;+\;\frac{N!}{(N+1)x^{N-1}}\right)x^{s+N-2}\,dx\;-\;\frac{N!}{s(N+1)}\\&&
\\& & \;+\;\int_1^\infty  \left(\frac{1}{e^x-T_{N-1}(x)}-\frac{N!}{x^N}\right)x^{s+N-2}\,dx.\end{eqnarray*}
Hence,
\begin{equation}\zeta_N(s)=\frac{1}{\Gamma(s+N-1)}\int_0^\infty \left(\frac{1}{e^x-T_{N-1}(x)}-\frac{N!}{x^N}\;+\;\frac{N!}{(N+1)x^{N-1}}\right)x^{s+N-2}dx.\label{eq 3.2}\end{equation} \end{rem}
From the above theorem and remark,  it may appear that hypergeometric zeta has an infinite number of poles since each application produces a pole on the right hand side of $\left(\ref{eq 3.2}\right)$; however, after $N$ repetitions the poles of $ \Gamma(s+N-1)$  on the left hand side begin to make their appearance, thereby canceling those on the right.  Hence, hypergeometric zeta has at most a finite number of poles.  We will have more to say about this in our second approach using contour integration (see Theorem 3.3).\\

The main advantage in using $\left(\ref{eq 3.1}\right)$ to analytically continue $\zeta_N(s)$  is that it reveals the behavior of $\zeta_N(s)$  near the pole $s=1$. This is the content of the next theorem.

\begin{tem} For $N\geq 1$, we have
\begin{equation}\lim_{s\rightarrow 1}\left[\zeta_N(s)-\frac{N!}{s-1}\right]= \log(N!)-N\frac{\Gamma'(N)}{\Gamma(N)}.\label{eq 3.3}\end{equation}
\end{tem}
\begin{proof} From $\left(\ref{eq 3.1}\right)$ we have
\begin{eqnarray*} \Gamma(s+N-1)\zeta_N(s)-\frac{N!}{s-1}= \int_0^1 \left(\frac{1}{e^x-T_{N-1}(x)}-\frac{N!}{x^N}\right)x^{s+N-2}\,dx \;+\;\int_1^\infty \frac{x^{s+N-2}}{e^x-T_{N-1}(x)}\,dx.\end{eqnarray*}
It follows from the Dominated Convergence Theorem that
\begin{eqnarray} \lim_{s\rightarrow 1}\left[\Gamma(s+N-1)\zeta_N(s)-\frac{N!}{s-1}\right]= \int_0^1 \left(\frac{1}{e^x-T_{N-1}(x)}-\frac{N!}{x^N}\right)x^{N-1}\,dx \;+\;\int_1^\infty \frac{x^{N-1}}{e^x-T_{N-1}(x)}\,dx.\label{eq 3.4}\end{eqnarray}
Now, 
\begin{eqnarray}
\int_0^1 \left(\frac{1}{e^x-T_{N-1}(x)}-\frac{N!}{x^N}\right)x^{N-1}\,dx&=&(N-1)!\left[\log\left(e^x-T_{N-1}(x)\right)-x-N\log x \right]\big|_0^1\nonumber\\
&=&(N-1)!\log\left[\frac{e^x-T_{N-1}(x)}{x^N e^x}\right]_0^1\label{eq 3.5}
\\ &=&(N-1)!\log\left[\frac{e-T_{N-1}(1)}{ e}\right]+(N-1)!\log(N!).\nonumber\end{eqnarray}
Also
\begin{eqnarray}
\int_1^\infty  \frac{x^{N-1}}{e^x-T_{N-1}(x)}\,dx&=&(N-1)!\left(\log \left[e^x-T_{N-1}(x) \right]-x \right)_1^\infty  \nonumber\\&=&(N-1)!\log\left[\frac{e^x-T_{N-1}(x)}{ e^x}\right]_1^\infty \label{eq 3.6}\\
&=&-(N-1)!\log\left[\frac{e-T_{N-1}(1)}{ e}\right].\nonumber
\end{eqnarray}
Using  $\left(\ref{eq 3.5}\right)$  and $\left(\ref{eq 3.6}\right)$ in $\left(\ref{eq 3.4}\right)$, we obtain
\begin{eqnarray*} \lim_{s\rightarrow 1}\left[\Gamma(s+N-1)\zeta_N(s)-\frac{N!}{s-1}\right]=(N-1)!\log(N!). \end{eqnarray*}
Hence,
\begin{eqnarray*}\lim_{s\rightarrow 1}\left[\zeta_N(s)-\frac{N!}{s-1}\right]&=& \lim_{s\rightarrow 1}\left[\frac{\left(\Gamma(s+N-1)\zeta_N(s)-\frac{N!}{s-1}\right)}{\Gamma(s+N-1)}\right ]-\lim_{s\rightarrow 1}\left[\frac{N\Gamma(s+N-1)-N!}{(s-1)\Gamma(s+N-1)}\right] \\ & &\\ &=& \frac{(N-1)!\log(N!)}{\Gamma(N)}-N\lim_{s\rightarrow 1}\left[\frac{1}{\Gamma(s+N-1)}\frac{\Gamma(s+N-1)-\Gamma(N)}{(s-1)}\right]\\& & \\ &=&\log(N!)-N\frac{\Gamma'(N)}{\Gamma(N)}.
\end{eqnarray*}
\end{proof}

\begin{rem} Observe that $(\ref{eq 3.3})$  yields the following classic result for $\zeta(s)$  (cf. \cite{T}):
\begin{equation}\lim_{s\rightarrow 1}\left[\zeta(s)-\frac{1}{s-1}\right]= -\frac{\Gamma'(1)}{\Gamma(1)}=\gamma\approx 0.577. \label{eq 3.7}\end{equation}\end{rem}

We now take a different approach and follow Riemann by using contour integration to develop the analytic continuation.  This will allow us to not only make precise our earlier statement about  $\zeta_N(s)$ having a finite number of poles but also to make explicit the role of the zeros of $e^x-T_{N-1}(x)=0$  in determining the values of hypergeometric zeta at negative integers. \\

\begin{figure} \centering
\includegraphics[width=.35\textwidth]{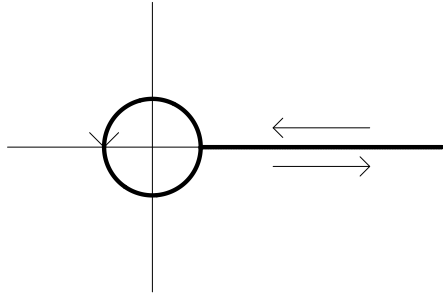}
\caption{Contour $\gamma$.}
\label{figure:contour-gamma}
\end{figure}

To this end consider the contour integral
\begin{eqnarray} I_N(s)=\frac{1}{2\pi i}\int_\gamma \left( e^w-T_{N-1}(w)\right)^{-1}(-w)^{s+N-1}\frac{dw}{w},\label{eq 3.8}\end{eqnarray}
where the contour  $\gamma$  is taken to be along the real axis from $\infty$  to $\delta>0$, then counterclockwise around the circle of radius $\delta$, and lastly along the real axis from $\delta$ to $\infty$ (cf. Figure \ref{figure:contour-gamma}).  Moreover, we let $-w$  have argument  $-\pi$ backwards along  $\infty$ to $\delta$  and argument $\pi$  when going to $\infty$.  Also, we choose the radius $\delta$  to be sufficiently small (depending on $N$) so that there are no roots of $e^w-T_{N-1}(w)=0$  inside the circle of radius $\delta$  besides the trivial root $z_0=0$.  This follows from the fact that $z_0=0$ is an isolated zero.  It is then clear from this assumption that  $I_N(s)$ must converge for all complex $s$ and therefore defines an entire function.\\

\begin{rem} \ 
\begin{enumerate}
\item[(a)] To be precise the contour  $\gamma$ should be taken as a limit of contours $\gamma_\epsilon$  as $\epsilon \rightarrow 0$, where the portions running along the $x$-axis are positioned at heights $\pm \epsilon$.  Moreover, the poles of the integrand in  $\left(\ref{eq 3.8}\right)$ cannot accumulate inside this strip due to the asymptotic exponential growth of the zeros of $e^w-T_{N-1}(w)=0$  (see Appendix I).
\item[(b)] Since we are most interested in the properties of $I_N(s)$ in the limiting case when $\delta \rightarrow 0$, we will also write $I_N(s)$ to denote $\lim_{\delta \rightarrow 0}I_N(s)$.  No confusion should arise from this abuse of notation.
\end{enumerate}
\end{rem}

We begin by evaluating $I_N(s)$  at integer values of $s$.  To this end, we decompose it as follows:\begin{eqnarray} I_N(s)&=&\frac{1}{2\pi i}\int_\infty^\delta \left( e^x-T_{N-1}(x)\right)^{-1}e^{(s+N-1)(\log x-\pi i)}\frac{dx}{x}\nonumber 
\\& & \;+\;\frac{1}{2\pi i}\int_{|w|=\delta} \left( e^w-T_{N-1}(w)\right)^{-1}(-w)^{s+N-1}\frac{dw}{w}\label{eq 3.9}\\
& &  \;+\;\frac{1}{2\pi i}\int_\delta^\infty\left( e^x-T_{N-1}(x)\right)^{-1}e^{(s+N-1)(\log x+i\pi)}\frac{dx}{x}.\nonumber\end{eqnarray}
Now, for integer  $s=n$, the two integrations along the real axis in $\left(\ref{eq 3.9}\right)$  cancel and we are left with just the middle integral around the circle of radius $\delta$:
\begin{eqnarray*} I_N(n)=\frac{1}{2\pi i}\int_{|w|=\delta} \left( e^w-T_{N-1}(w)\right)^{-1}(-w)^{n+N-1}\frac{dw}{w}.\end{eqnarray*}
Since the expression  $w^N\left(e^w-T_{N-1}(w)\right)^{-1}$ inside the integrand has a removable singularity at the origin, it follows by Cauchy's Theorem that for integers $n>1$, $$I_N(n)=0.$$ For integers $n \leq 1$, we consider the power series expansion
\begin{eqnarray} \frac{w^N/N!}{ e^w-T_{N-1}(w)}=\sum_{m=0}^\infty\frac{B_{N,m}}{m!}w^m.\label{eq 3.10}\end{eqnarray}
 It now follows from the Residue Theorem that
\begin{eqnarray} I_N(n)&=&\frac{1}{2\pi i}\int_{|w|=\delta} \left( e^w-T_{N-1}(w)\right)^{-1}(-w)^{n+N-1}\frac{dw}{w}\nonumber\\
&=&(-1)^{n+N-1}\frac{N!}{2\pi i}\int_{|w|=\delta} \left(\sum_{m=0}^\infty\frac{B_{N,m}}{m!}w^m\right)\frac{dw}{w^{2-n}}\label{eq 3.11}\\
&=&\frac{(-1)^{n+N-1}N!B_{N,1-n}}{(1-n)!}.\nonumber\end{eqnarray}

We now express $\zeta_N(s)$ in terms of $I_N(s)$.  For $\Re(s)=\sigma >1$, the middle integral in $\left(\ref{eq 3.9}\right)$ goes to zero as $\delta \rightarrow 0$.  It follows that
\begin{eqnarray*} I_N(s)&=&\left(\frac{e^{\pi i(s+N-1)}-e^{-\pi i(s+N-1)} }{2\pi i}\right)\int_0^\infty \left( e^x-T_{N-1}(x)\right)^{-1}x^{s+N-2}dx\\ & & \\
& =&\frac{\sin\left[\pi(s+N-1)\right]}{\pi }\Gamma(S+N-1)\zeta_N(s).   \end{eqnarray*}
Now, by using the functional equation for the gamma function:$$\Gamma(1-(s+N-1))\Gamma(s+N-1)=\frac{\pi}{\sin[\pi(s+N-1)]}$$ we obtain 
\begin{eqnarray} \zeta_N(s)=\Gamma(1-(s+N-1))I_N(s). \label{eq 3.12}\end{eqnarray}

\begin{rem} Equation $\left(\ref{eq 3.12}\right)$ implies that the zeros of $I_N(s)$  at positive integers $n>1$  are simple since we know from Theorem 2.1  that  $\zeta_N(n)>1$ for $n>1$. \end{rem}

 Here is another consequence of $\left(\ref{eq 3.12}\right)$, which we state as 
\begin{tem} $\zeta_N(s)$  is analytic on the entire complex plane except for simple poles at $\{2-N,\, 3-N,\; \cdots, \;1\}$   whose residues are
\begin{eqnarray} {\mbox Res}\left(\zeta_N(s),s=n\right)=(2-n)\left( \begin{array}{cc} N  \\ 2-n \end{array} \right)B_{N,1-n}\;\;\;\;\;\;\;(2-N\leq n\leq 1).\label{eq 3.13} \end{eqnarray}
Further more, for negative integers  $n$ less than $2-N$, we have \begin{eqnarray} \zeta_N(n)=(-1)^{-n-N+1}\left( \begin{array}{cc} 1-n  \\ N \end{array} \right)^{-1}B_{N,1-n}. \label{eq 3.14}\end{eqnarray}\end{tem}
\begin{proof} Since   $\Gamma(1-(s+N-1))$ has only simple poles at $s=2-N,\, 3-N,\; \cdots $, and $I_N(s)$ has simple zeros at $s=2,\;3,\, \cdots \;$, it follows from  $\left(\ref{eq 3.12}\right)$ that $\zeta_N(s)$  is analytic on the whole plane except for simple poles at  $s=n$, $2-N \leq n \leq 1$.  Recalling the fact that the residue of  $\Gamma(s)$ at negative integer $n$ is $(-1)^n/|n|!$, it follows from $\left(\ref{eq 3.11}\right)$ that the residue of $\zeta_N(s)$  at the same pole is:
\begin{eqnarray*} {\mbox Res}\left(\zeta_N(s),s=n\right)&=&\lim_{s\rightarrow n}(s-n)\zeta_N(s)\;=\;\lim_{s\rightarrow n}\left[(s-n)\Gamma(1-(s+N-1))I_N(s)\right]\\& & \\
&=&-\frac{(-1)^{2-N-n}}{(2-N-n)!}I_N(n)\;=\;-\frac{(-1)^{2-N-n}}{(2-N-n)!}\frac{(-1)^{n+N-1}N!B_{N,1-n}}{(1-n)!}\\& & \\
&=&(2-n)\left( \begin{array}{cc} N  \\ 2-n \end{array} \right)B_{N,1-n},  \end{eqnarray*}
which proves $\left(\ref{eq 3.13}\right)$.  For $n<2-N$, $\left(\ref{eq 3.11}\right)$, $\left(\ref{eq 3.12}\right)$, and the fact that $\Gamma(1-(n+N-1))=(1-N-n)!$ imply
\begin{eqnarray*} \zeta_N(n)&=&\Gamma(1-(n+N-1)I_N(n)=\frac{(-1)^{n+N-1}N!(1-N-n)!B_{N,1-n}}{(1-n)!}\\&=&(-1)^{-n-N+1}\left( \begin{array}{cc} 1-n  \\ N \end{array} \right)^{-1}B_{N,1-n}, \end{eqnarray*} which is $\left(\ref{eq 3.14}\right).$ This completes the proof the theorem.
\end{proof}

\begin{rem} \ 
\begin{enumerate}
\item[(a)] We note that the coefficients $B_{N,n}$ defined by $\left(\ref{eq 3.10}\right)$ generalize the Bernoulli numbers $B_n$, which arise when $N = 1$.  For $N = 2$, the coefficients $B_{2,n}$  have been studied extensively by Howard \cite{H1}-\cite{H3}, who referred to them as $A_n$.  We will use some of Howard's results in the next section to obtain a functional inequality (as opposed to a functional equation) involving $\zeta_2(s)$.   For $N$ in general, we note that the coefficients $B_{N,n}$ can be found recursively by the relation
\begin{eqnarray*}  B_{N,0}=1,\;\;\;\;\;\;\sum_{m=0}^n\frac{n!B_{N,m}}{(N+n-m)!m!}=0 \;\;(n\geq 1).  \end{eqnarray*} Or equivalently, \begin{eqnarray*}  B_{N,0}=1,\;\;\;\;\;\;B_{N,n}=-\sum_{m=0}^{n-1}\frac{n!B_{N,m}}{(N+n-m)!m!} \;\;(n\geq 1).  \end{eqnarray*}
Here are the first few values of $B_{N,n}$:
 \begin{eqnarray*}\begin{array}{ll} B_{N,0}=1,& B_{N,1}=- \frac{1}{N+1}, \\ B_{N,2}= \frac{2}{(N+1)^2(N+2)}, & B_{N,3}= \frac{6}{(N+1)^3(N+2)(N+3)}.\end{array}\end{eqnarray*}
\item[(b)] It follows that the residues of $\zeta_N(s)$  can be found similarly by recursion.  For example:
\begin{eqnarray*} {\rm{Res}}\left(\zeta_N(s),s=n\right)=\left\{ \begin{array}{cc} N&n=1, \\& \\ -\frac{N(N-1)}{N+1}&n=0, \\&\\\frac{N(N-1)(N-2)}{(N+1)^2(N+2)}&n=-1,  \\&\\\frac{N(N-1)^2(N-2)(N-3)}{(N+1)^3(N+2)(N+3)}&n=-2.\end{array}\right. \end{eqnarray*}
\end{enumerate}
\end{rem}

We end this section with the following result: 
\begin{tem}
\begin{equation}I_N'(1)=(-1)^N(N-1)!\log(N!).\label{eq 3.15}\end{equation}
\end{tem}  
  \begin{proof}  Following Edwards in \cite{E}, we rewrite $I_N(s)$ as follows:
\begin{eqnarray*} I_N(s)&=&\frac{1}{2\pi i}\int_\infty^\delta \frac{\left( x e^{-\pi i}\right)^{(s+N-1)}}{e^x-T_{N-1}(x)}\frac{dx}{x} \;+\;\frac{1}{2\pi i}\int_{|w|=\delta}  \frac{(-w)^{s+N-1}}{ e^w-T_{N-1}(w)}\frac{dw}{w} \;+\;\frac{1}{2\pi i}\int_\delta^\infty\frac{\left(x  e^{\pi i}\right)^{(s+N-1)}}{e^x-T_{N-1}(x)}\frac{dx}{x}.\end{eqnarray*}
\noindent It follows that
\begin{eqnarray*} I_N'(s)&=&\frac{1}{2\pi i}\int_\infty^\delta \frac{\left(x e^{-\pi i}\right)^{(s+N-1)}(\log x -i\pi)}{e^x-T_{N-1}(x)}\frac{dx}{x} \;+\;\frac{1}{2\pi i}\int_{|w|=\delta}  \frac{(-w)^{s+N-1}(\log \delta+i\theta -i\pi)}{ e^w-T_{N-1}(w)}\frac{dw}{w}\\& & \\& &  \;+\;\frac{1}{2\pi i}\int_\delta^\infty\frac{\left( x e^{\pi i}\right)^{(s+N-1)}(\log x + i\pi)}{e^x-T_{N-1}(x)}\frac{dx}{x}.\end{eqnarray*}
Therefore, 
\begin{eqnarray*} I_N'(1)&=&\frac{1}{2\pi i}\int_\infty^\delta \frac{\left(x e^{-\pi i}\right)^N(\log x -i\pi)}{e^x-T_{N-1}(x)}\frac{dx}{x} \;+\;\frac{1}{2\pi i}\int_{|w|=\delta}  \frac{(-w)^N(\log \delta+i\theta -i\pi)}{ e^w-T_{N-1}(w)}\frac{dw}{w}\\& & \\& &  \;+\;\frac{1}{2\pi i}\int_\delta^\infty\frac{\left( x e^{\pi i}\right)^N(\log x + i\pi)}{e^x-T_{N-1}(x)}\frac{dx}{x}\\& & \\
&=&(-1)^N \int_\delta^\infty \frac{  x^{N-1}}{e^x-T_{N-1}(x)}dx \;+\;(-1)^N\frac{\log \delta}{2\pi i}\int_{|w|=\delta}  \frac{(-w)^N}{ e^w-T_{N-1}(w)}\frac{dw}{w}\\& & \\& &  \;-\;\frac{1}{2\pi i}\int_{-\pi}^{\pi}\frac{(-w)^N}{ e^w-T_{N-1}(w)}\phi d\phi.\end{eqnarray*}
Now, the third integral on the right hand side approaches zero as $\delta \rightarrow 0$.  As for the other two integrals, observe that the first is the same as $\left(\ref{eq 3.6}\right)$, except for the lower limit of integration.  We proceed as before and observe that
\begin{eqnarray*}
\int_\delta^\infty  \frac{x^{N-1}}{e^x-T_{N-1}(x)}\,dx&=&(N-1)!\log\left[\frac{e^x-T_{N-1}(x)}{ e^x}\right]_\delta^\infty \\
&=&-(N-1)!\log\left(\frac{e^\delta-T_{N-1}(\delta)}{e^\delta}\right)\\&=&-(N-1)!\log\left(\frac{\delta^N}{N!}+\frac{\delta^{N+1}}{(N+1)!}+\cdots\right)\\
&=&-N!\log \delta -(N-1)!\log\left(\frac{1}{N!}+\frac{\delta}{(N+1)!}+\cdots\right).\end{eqnarray*}
The second integral we have already encountered and can be evaluated using residue theory:
\begin{eqnarray*}\frac{\log \delta}{2\pi i}\int_{|w|=\delta}  \frac{(-w)^N}{ e^w-T_{N-1}(w)}\frac{dw}{w}=N!\log \delta.\end{eqnarray*}
Therefore, in the limit as  $\delta \rightarrow 0$, we obtain $\left(\ref{eq 3.15}\right)$.
\end{proof}

\begin{rem} Observe that Theorem 3.4 by itself does not yield the classical result
	 \begin{eqnarray}\frac{\zeta'(0)}{\zeta(0)}=2\pi.\label{eq 3.16} \end{eqnarray}
As Edwards demonstrates in \cite{E} the proof of $\left(\ref{eq 3.16}\right)$ also relies on the functional equation for $\zeta(s)$.  Therefore, it is unclear how $\left(\ref{eq 3.16}\right)$ generalizes to an analogous formula for $\zeta_N'(s)/\zeta_N(s)$ at suitable negative integer values of $s$ since no functional equation is known for $\zeta_N(s)$ when $ N  > 1$ . 
\end{rem}

\section{ FUNCTIONAL INEQUALITY}

\setcounter{equation}{0}   

\begin{figure} \centering
\includegraphics[width=.35\textwidth]{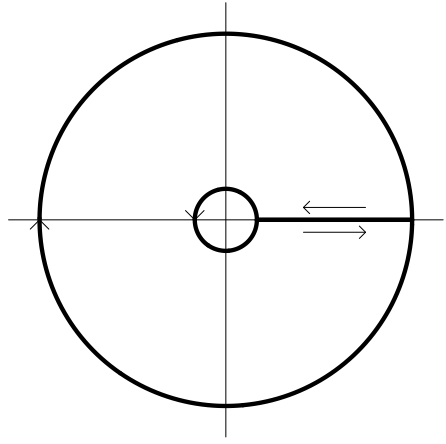}
\caption{Contour $\gamma_M$.}
\label{figure:contour-gamma-m}
\end{figure}

 In the present section, we discuss a `functional inequality' satisfied by $\zeta_N(s)$. Let $\gamma_M $  be the annulus-shaped contour consisting of two concentric circles centered at the origin, the outer circle having radius $(2M+1)\pi$ and the inner circle having radius $\delta <\pi$ (cf. Figure \ref{figure:contour-gamma-m}).  The outer circle is traversed clockwise, the inner circle counterclockwise and the radial segment along the positive real axis is traversed in both directions.  Then define
\begin{eqnarray} I_{\gamma_M}(s)=\frac{1}{2\pi i}\int_{\gamma_M} \frac{(-z)^{s+N-1}}{e^z-T_{N-1}(z)}\frac{dz}{z}.\label{eq 4.1}\end{eqnarray}

We claim that  $I_{\gamma_M}(s)$ converges to $I_N(s)$ as $M\rightarrow\infty$ for $\Re(s)<0$. To prove this, observe that the portion of $I_{\gamma_M}(s)$ around the outer circle tends to zero as $M\rightarrow\infty$   on the same domain.  This is because on the outer circle defined by $|z|=(2M+1)\pi$ we have that  $\left |z^{N-1}/(e^z-T_{N-1}(z)\right|$ is bounded independently of $M$ and $\left |(-z)^s/z\right|<|z|^{\Re(s)-1}$.  Therefore,
\begin{eqnarray} I_N(s)=\lim_{M\rightarrow\infty}I_{\gamma_M}(s).\label{eq 4.2}\end{eqnarray}

On the other hand, we have by residue theory
\begin{eqnarray} I_{\gamma_M}(s)=-\sum_{k=1}^K\left[\mbox{Res}\left(\frac{(-z)^{s+N-2}}{e^z-T_{N-1}(z)},z=z_k\right)\;+\;\mbox{Res}\left(\frac{(-z)^{s+N-2}}{e^z-T_{N-1}(z)},z=\bar{z}_k\right)\right]. \label{eq 4.3}\end{eqnarray}
Here, $z_k=r_k e^{i\theta_k}$ and $\bar{z}_k=r_k e^{-i\theta_k}$ are the complex conjugate roots of $e^z-T_{N-1}(z)=0$ and $K=K_M$ is the number of roots inside $\gamma_M$ in the upper-half plane. Clearly $z_k$  depends on $N$. We will make this assumption throughout  and use the same notation $z_k$ instead of the more cumbersome notation $z_k(N)$.   Moreover, we arrange the roots in ascending order so that $|z_1|<|z_2|<|z_3|<\cdots$, since none of the roots can have the same length (see Appendix I).  Now, to evaluate the residues, we call upon Cauchy's Integral Formula:
   
\begin{eqnarray*}\mbox{Res}\left(\frac{(-z)^{s+N-2}}{e^z-T_{N-1}(z)},z=z_k\right) = \int_{C_{k}}\left[\frac{(-z)^{s+N-2(z-z_k)}}{e^z-T_{N-1}(z)}\right]\frac{dz}{z-z_k} = (-z_k)^{s+N-2}\lim_{z\rightarrow z_k}\frac{z-z_k}{e^z-T_{N-1}(z)}.\end{eqnarray*} Here, $C_k$ is any sufficiently small contour enclosing only one root $z_k$ of $e^z-T_{N-1}(z)=0$.  But then
\begin{eqnarray*}\lim_{z\rightarrow z_k}\frac{z-z_k}{e^z-T_{N-1}(z)}=\frac{1}{e^{z_k}-T_{N-2}(z_k)}=\frac{(N-1)!}{z_k^{N-1}}.\end{eqnarray*}
It follows that\begin{eqnarray*}\mbox{Res}\left(\frac{(-z)^{s+N-2}}{e^z-T_{N-1}(z)},z=z_k\right) =  (-1)^{N-1}(N-1)!(-z_k)^{s-1}.\end{eqnarray*}
Therefore, 
\begin{eqnarray}I_{\gamma_M}(s)&=&(-1)^{N-1}(N-1)!\sum_{k=1}^K\left[ (-z_k)^{s-1}+(-\bar{z}_k)^{s-1}\right]\nonumber\\&=&  2(-1)^{N-1}(N-1)!\sum_{k=1}^Kr_k^{s-1}\cos\left[(s-1)(\pi-\theta_k)\right].\label{eq 4.4}\end{eqnarray}Since $K\rightarrow \infty$ as $M\rightarrow \infty$, we have by $\left(\ref{eq 4.2}\right)$ and $\left(\ref{eq 4.4}\right)$, \begin{eqnarray}I_N(s)&=&\lim_{M\rightarrow \infty}
I_{\gamma_M}(s)\nonumber\\&=&2(-1)^{N-1}(N-1)!\sum_{k=1}^\infty r_k^{s-1}\cos\left[(s-1)(\pi-\theta_k)\right].\label{eq 4.5}\end{eqnarray} Combining    $\left(\ref{eq 3.12}\right)$   and $\left(\ref{eq 4.5}\right)$  we have proved

\begin{tem}For $\Re(s)<0$,\begin{eqnarray} \zeta_N(s)=2(-1)^{N-1}(N-1)!\Gamma(1-(s+N-1))\sum_{k=1}^\infty r_k^{s-1}\cos\left[(s-1)(\pi-\theta_k)\right].\label{eq 4.6}\end{eqnarray} 
\end{tem}

\begin{rem} \
\begin{enumerate}
\item[(a)] Observe that for $N=1$ equation  $\left(\ref{eq 4.6}\right)$ reduces to the classical functional equation since in this case we have $z_k=2\pi k i$, and therefore $r_k=2\pi k$ and $\theta_k=\pi/2$:
\begin{eqnarray} \zeta(s)= 2(2\pi)^{s-1}\sin\left(\frac{\pi}{2}s\right)\Gamma(1-s)\zeta(s).\label{eq 4.7}\end{eqnarray}
\item[(b)] The first 10 nonzero roots $\{z_k\}$ of $e^z-T_{N-1}(z)=0$ are listed in Appendix II for the cases $N=2$ and $N=3$. 
\end{enumerate}
\end{rem}
Next we establish a connection between $\zeta_2(s)$ and the classical zeta function. More specifically, we prove $\left(\ref{eq 1.6}\right)$, which we restate as 
\begin{tem}For $\Re(s)<0$, we have \begin{equation}\left| \zeta_2(s)\right|<2(2\pi)^{\Re(s)}\left|\Gamma(-s)\right | e^{\Im(s)(\pi-\theta_1)}\zeta(1-\Re(s)) . \label{eq 4.8} \end{equation} 
\end{tem}
\begin{proof} The argument essentially rests on bounds obtained by Howard on the zeros of  $e^z-T_{N-1}(z)=0$.  Since he provides few details of the proof in \cite{H3}, we give a full proof of it in Appendix I (Lemma 5.3).  In particular, we will establish  that there are infinitely many zeros and that all of them are simple.  Moreover, for each positive integer $k$, there exists precisely one zero $z=x+i y=re^{i\theta}$ whose imaginary part is bounded by
\begin{equation}(2k+1/4)\pi\;<\;y\;<\;(2k+1/2)\pi.\label{eq 4.9}\end{equation} It follows that the zeros satisfying $\left(\ref{eq 4.9}\right)$, their conjugates, and $z = 0$ exhaust all the zeros of $e^z-T_{N-1}(z)=0$.  We then order the zeros $\{z_k\}$ in the upper-half plane so that $0=|z_0|=<|z_1|<|z_2|<\cdots $   (see Appendix I).\\
Now, using the fact that the angles $\{\theta_k\}$  are monotonically increasing ( to $\pi/2$), we have \begin{eqnarray*}|\cos\left[(s-1)(\pi-\theta_k)\right]|\leq e^{|\Im(s)(\pi-\theta_k)|}<e^{|\Im(s)(\pi-\theta_1)|}.\end{eqnarray*}
Therefore, by $\left(\ref{eq 4.6}\right)$ and $\left(\ref{eq 4.9}\right)$,  for $\Re(s)<0$, the following bound is achieved:
\begin{eqnarray*} |\zeta_2(s)|&=&\left|2\Gamma(-s)\sum_{k=1}^\infty r_k^{s-1}\cos\left[(s-1)(\pi-\theta_k)\right]\right|\\&\leq & 2e^{|\Im(s)(\pi-\theta_1)|}|\Gamma(-s)|\sum_{k=1}^\infty \frac{1}{r_k^{1-\Re(s)}}\\&\leq&2e^{|\Im(s)(\pi-\theta_1)|}|\Gamma(-s)|\sum_{k=1}^\infty \frac{1}{y_k^{1-\Re(s)}}\\&< &2e^{|\Im(s)(\pi-\theta_1)|}|\Gamma(-s)|\sum_{k=1}^\infty \frac{1}{(2\pi k)^{1-\Re(s)}}\\&=&2(2\pi)^{\Re(s)-1}e^{|\Im(s)(\pi-\theta_1)|}|\Gamma(-s)|\sum_{k=1}^\infty \frac{1}{ k^{1-\Re(s)}}\\&=&2(2\pi)^{\Re(s)-1}e^{|\Im(s)(\pi-\theta_1)|}|\Gamma(-s)|\zeta(1-\Re(s)).\end{eqnarray*}
This completes the proof.
\end{proof}

Since $\zeta(1-\Re(s))\leq \zeta_2(1-\Re(s))$  from $\left(\ref{eq 2.10}\right)$, we thus obtain as a corollary the following functional inequality for $\zeta_2(s)$:
\begin{cor}For  $\Re(s)<0$, we have
 \begin{equation}  \left| \zeta_2(s)\right|<2(2\pi)^{\Re(s)}\left|\Gamma(-s)\right | e^{\Im(s)(\pi-\theta_1)}\zeta_2(1-\Re(s)). \label{eq 4.10} \end{equation}\end{cor}

\begin{rem} \
\begin{enumerate}
\item[(a)] Compare $\left(\ref{eq 4.10}\right)$ with the functional equation $\left(\ref{eq 4.7}\right)$ of the Riemann zeta function.
\item[(b)] We can improve on $\left(\ref{eq 4.8}\right)$ using the bound from $\left(\ref{eq 4.9}\right)$, namely $(2k+1/4)\pi<b_k$:
\begin{eqnarray*} |\zeta_2(s)|&
<&2e^{|\Im(s)(\pi-\theta_1)|}|\Gamma(-s)|\sum_{k=1}^\infty \frac{1}{y_k^{1-\Re(s)}}\\&
<&2e^{|\Im(s)(\pi-\theta_1)|}|\Gamma(-s)|\sum_{k=1}^\infty \frac{1}{((2k+1/4)\pi)^{1-\Re(s)}}\\&<&2(2\pi)^{\Re(s)-1}e^{|\Im(s)(\pi-\theta_1)|}|\Gamma(-s)|\sum_{k=1}^\infty \frac{1}{ (k+1/8)^{1-\Re(s)}}\\&
<&2(2\pi)^{\Re(s)-1}e^{|\Im(s)(\pi-\theta_1)|}|\Gamma(-s)|\zeta(1-\Re(s),1/8).\end{eqnarray*} Here, $\zeta(s,a)$ is the Hurwitz zeta function defined by
$$\zeta(s,a)=\sum_{n=1}^{\infty}\frac{1}{(n+a)^s}.$$
\end{enumerate}
\end{rem}

For an asymptotically tighter bound where $s$ is sufficiently large, we make use of the inequality \begin{equation}r_k\geq  m r_1, \;\mbox{where}\; k=2m \;\mbox{or} \;k=2m-1.\label{eq 4.11}\end{equation}
To prove $\left(\ref{eq 4.11}\right)$, we observe that $r_1<7.8$ and $r_k>y_k>(2k+1/4)\pi$. Now if $k=2m$, then $r_k>(4m+1/4)\pi>4m \pi >m r_1 .$ So suppose  $k=2m-1$.  For $m=1$, the assertion is obvious. If  $m\geq 2$, then $m>7\pi/(4(4\pi-r_1))$ and hence    $r_k>(4m\pi -7\pi/4)>m r_1.$ This completes the proof of $\left(\ref{eq 4.11}\right)$. \\
It now follows from $\left(\ref{eq 4.11}\right)$ that
\begin{eqnarray*}\sum_{k=1}^\infty\frac{1}{r_k^{1-\Re(s)}}&=&\frac{1}{r_1^{1-\Re(s)}}\sum_{k=1}^\infty\left(\frac{r_1}{r_k}\right)^{1-\Re(s)}\nonumber\\&\leq &\frac{1}{r_1^{1-\Re(s)}}\left[\sum_{m=1}^\infty\left(\frac{r_1}{r_{2m}}\right)^{1-\Re(s)}\;+\;\sum_{m=1}^\infty\left(\frac{r_1}{r_{2m-1}}\right)^{1-\Re(s)}\right]\nonumber\\& & \\&< &\frac{1}{r_1^{1-\Re(s)}}\left[\sum_{m=1}^\infty\left(\frac{1}{m}\right)^{1-\Re(s)}\;+\;\sum_{m=1}^\infty\left(\frac{1}{m}\right)^{1-\Re(s)}\right]\\
& = & \frac{2\zeta(1-\Re(s))}{r_1^{1-\Re(s)}}. \end{eqnarray*}
This produces the following bound on $\zeta_2(s)$:\begin{tem} \begin{equation}\left| \zeta_2(s)\right|<4r_1^{\Re(s)-1}\left|\Gamma(-s)\right | e^{\Im(s) (\theta_1-\pi)}\zeta(1-\Re(s)). \label{eq 4.12} \end{equation} 
\end{tem}

 We now assume that $s=-(n-1)$  is a negative integer less than $1-N$.  It then follows from $\left(\ref{eq 4.6}\right)$ and  $\left(\ref{eq 3.14}\right)$ that \begin{eqnarray*}B_{N,n}=(-1)^{N-1}\frac{2n!}{N}\sum_{k=1}^\infty r_k^{-n}\cos[n\theta_k]. \end{eqnarray*}
Observe that when $N = 2$ we obtain Howard's result in \cite{H3}: \begin{eqnarray*}B_{2,n}=-n!\sum_{k=1}^\infty r_k^{-n}\cos[n\theta_k].  \end{eqnarray*}Moreover, \begin{eqnarray*}\left|B_{N,n}\right|<\frac{2n!}{N(2\pi)^n} \zeta(n).  \end{eqnarray*} 

Since $\zeta(n)\leq\zeta(2)=\pi^2/6$, this establishes the following bound on the generalized Bernoulli numbers:
\begin{tem}(Howard \cite{H3}) For positive integers $n>N$, \begin{eqnarray}\left|B_{N,n}\right|<\frac{2n!}{N(2\pi)^n} \frac{\pi^2}{6}. \label{eq 4.13}\end{eqnarray}\end{tem}

\begin{rem} In the case when $N = 2$, $\left(\ref{eq 4.12}\right)$ and $\left(\ref{eq 3.14}\right)$ can be combined to improve on the bound $\left(\ref{eq 4.13}\right)$:\begin{eqnarray}\left|B_{2,n}\right|<\frac{2n!}{r_1^n}. \label{eq 4.14}  \end{eqnarray} Since $r_1\approx7.748>7$, this proves Howard's conjecture as stated in \cite{H3}: \begin{eqnarray*}\left|B_{2,n}\right|<\frac{n!}{7^n}.\end{eqnarray*}
Since the radius of convergence of the power series $\left(\ref{eq 3.10}\right)$ is  $|z_1|=r_1$, we note that $\left(\ref{eq 4.14}\right)$  is sharp and that  $\{B_{2,n}\}$ is an unbounded sequence.\end{rem}

\begin{rem}[Concluding Remarks] \
\begin{enumerate}
\item[(a)] In a work in progress we have proved that $\zeta_2(s)$ is zero-free on a suitable left half-plane. While there are ``trivial" zeros on the negative real axis, we have not yet been able to find any nonreal zeros of $\zeta_2(s)$.
\item[(b)] The theory of hypergeometric zeta functions can of course be extended to continuous values of the parameter $N$ via the definition
	  \begin{equation}\zeta(N;s)=\frac{\Gamma(N)}{\Gamma(s+N-1)}\int_0^\infty \frac{x^{s-1}}{_1F_1(1,N;x)-1}dx. \label {eq 4.15}\end{equation}
It is a straightforward exercise to verify that for integer values of $N$, $\left(\ref{eq 4.15}\right)$ reduces to our original definition of hypergeometric zeta functions given by $\left(\ref{eq 2.1}\right)$.  Observe that $\left(\ref{eq 4.15}\right)$ naturally leads to a continuous version of generalized Bernoulli numbers, a topic that has already investigated by K. Dilcher in \cite{D}.  We take up the theory of hypergeometric zeta functions based on $\left(\ref{eq 4.15}\right)$  in an upcoming paper.
\end{enumerate}
\end{rem}

\section{APPENDIX I}

\setcounter{equation}{0} 
In this appendix, we will investigate the roots of  \begin{equation}e^z-T_N(z)=0,\label{eq 5.1}\end{equation} where $T_N(z)=\sum_{k=0}^N\frac{z^k}{k!}$ and $N$ is a fixed positive integer. We shall prove our results in a sequence of lemmas. Lemma 5.1 gives an asymptotic approximation of the roots.  From its asymptotic formula one can conclude that these zeros can be arranged in an increasing order of magnitude.  Lemma 5.2 guarantees the existence of an infinite number of simple roots.   In Lemmas 5.3 and 5.4, we specialize to the cases $N=2$ and $N=3$.  In particular in Lemma 5.3, we will give a proof of the result of Howard \cite{H1} that we have used in Section 4.  Note that the functional equation of $\zeta_2(s)$, or the lack of it, depends on a detailed knowledge of these roots.  

\begin{lem}  Let $R>N$ be a positive real number and $\epsilon =N/R$. Define 
\begin{eqnarray*}A=\frac{1}{N!}\left(\frac{1-\epsilon^{N+1}}{1-\epsilon}\right), \;\;A_1=\sqrt{A^{-2/N}-\frac{1}{R}}, \;\;B= \frac{1}{N!}\left(2-\frac{1-\epsilon^{N+1}}{1-\epsilon}\right),\;\;\mbox{and} \;\; B_1=B^{-\frac{1}{N}}.\end{eqnarray*}Let $z=x+i y=re^{i \theta}$ be a root of $\left(\ref{eq 5.1}\right)$ with $y>0$. Then for $r>R$, we have \\

\noindent (i) $\;\;$ $B|z|^N\leq |T_N(z)|\leq A|z|^N$.\\
\noindent (ii) $\;\;$If $x>N \log R$, then \begin{eqnarray}A_1e^{\frac{x}{N}}\leq y\leq B_1e^{\frac{x}{N}}.\label{eq 5.2}\end{eqnarray}
\noindent (iii) \begin{eqnarray} y=2q\pi +N\theta +\delta_R,\label{eq 5.3}\end{eqnarray}where $q$ is an integer and $\delta$ is a real number such that $\lim_{R\rightarrow \infty}\delta_R=0$.
\end{lem}

\begin{rem}  When $R$ is sufficiently large, we note that $A\approx B\approx 1/N!$ and hence $A_1\approx B_!\approx (N!)^{1/N}$.  We deduce from $\left(\ref{eq 5.2}\right)$ that $y\approx (N!)^{1/N}e^{x/N}$.  This in turn implies  that $\theta \approx \pi/2$. Thus we have the asymptotic approximation of the roots of $\left(\ref{eq 5.1}\right)$:
\begin{eqnarray} x\approx N\log\left(2q\pi +N \frac{\pi}{2}-\log (N!)\right)\;\;\; \mbox{and}\; \;\; y\approx (N!)^{1/N}e^{x/N}. \label{eq 5.4}\end{eqnarray}\end{rem}

\begin{proof}[Proof of Lemma 5.1] For $|z|>R$, we apply the triangle inequality to obtain the upper bound 
\begin{eqnarray*}
\left|\frac{T_N(z)}{z^N}\right|& = &\left |\frac{1}{z^N}+\frac{1}{z^{N-1}}+\frac{1}{2!z^{N-2}}+\cdots +\frac{1}{(N-1)!z}+\frac{1}{N!}\right| \\& \leq & \frac{1}{|z|^N}+\frac{1}{|z|^{N-1}}+\frac{1}{2!|z|^{N-2}}+\cdots +\frac{1}{(N-1)!|z|}+\frac{1}{N!} \\& < & \frac{1}{N!}\left[\frac{N!}{R^N}+\frac{N!}{R^{N-1}}+\frac{N!}{2!R^{N-2}}+\cdots +\frac{N!}{(N-1)!R}+1\right]
\\& < &\frac{1}{N!}\left[\left(\frac{N}{R}\right)^{N}+\left(\frac{N}{R}\right)^{N-1}+\left(\frac{N}{R}\right)^{N-2}+\cdots +\frac{N}{R}+1\right]\\&=&A.\end{eqnarray*}
Similarly, using the triangle inequality in reverse, we obtain the lower bound 
\begin{eqnarray*}
\left|\frac{T_N(z)}{z^N}\right|& = &\left |\frac{1}{z^N}+\frac{1}{z^{N-1}}+\frac{1}{2!z^{N-2}}+\cdots +\frac{1}{(N-1)!z}+\frac{1}{N!}\right| \\& \geq & \frac{1}{N!}-\left[\frac{1}{|z|^N}+\frac{1}{|z|^{N-1}}+\frac{1}{2!|z|^{N-2}}+\cdots +\frac{1}{(N-1)!|z|}\right] \\& > & \frac{1}{N!}\left\{1-\left[ \left(\frac{N}{R}\right)^{N}+\left(\frac{N}{R}\right)^{N-1}+\left(\frac{N}{R}\right)^{N-2}+\cdots +\frac{N}{R}   \right]\right\}\\& = &\frac{1}{N!}\left\{2-\left[ \left(\frac{N}{R}\right)^{N}+\left(\frac{N}{R}\right)^{N-1}+\left(\frac{N}{R}\right)^{N-2}+\cdots +\frac{N}{R}+1   \right]\right\}\\&=&B.\end{eqnarray*}
This proves (i). To prove (ii), we note that since $e^z-T_N(z)=0$ and  $\left |e^z\right |=e^x$, (i) yields  \begin{eqnarray}B|z|^N\leq e^x \leq A|z|^N. \label{eq 5.5}\end{eqnarray}
Taking the $N^{th}$  root and squaring  $\left(\ref{eq 5.5}\right)$, we get $B^{2/N}|z|^2\leq e^{2x}{N}\leq B^{2/N}|z|^2$. We now solve these inequalities for $y$ to obtain \begin{eqnarray}e^{\frac{x}{N}}\sqrt{A^{\frac{-2}{N}}-x^2e^{\frac{-2x}{N}}} \leq y \leq e^{\frac{x}{N}}\sqrt{B^{\frac{-2}{N}}-x^2e^{\frac{-2x}{N}}}. \label{eq 5.6}\end{eqnarray}Since $x^2e^{-2x/N}$  is always positive the second inequality in (ii) follows from $\left(\ref{eq 5.6}\right)$. Note also that $x^2e^{-2x/N}$   is decreasing on  $[N,\infty)$ and thus for $x>N \log R$, we have \begin{eqnarray}x^2e^{\frac{-2x}{N}}< \left(N\log R\right)^2e^{\frac{-2 N\log R}{N}}=\frac{\left(N\log R\right)^2}{R^2}<\frac{1}{R}. \label{eq 5.7}\end{eqnarray} The first inequality of (ii) now follows from $\left(\ref{eq 5.6}\right)$ and $\left(\ref{eq 5.7}\right)$,  thereby establishing both inequalities. To prove (iii),  we observe that since $e^z=T_N(z)$ and $\arg\left(e^z\right)=\arg\left(e^{x+iy}\right)=y$ we have \begin{eqnarray}\arg\left(T_N(z)\right)=y. \label{eq 5.8}\end{eqnarray}On the other hand, for large $R$, we have \begin{eqnarray*} \arg\left(\frac{T_N(z)}{z^N}\right) = \arg \left (\frac{1}{z}\left[\frac{1}{z^{N-1}}+\frac{1}{z^{N-2}}+\frac{1}{2!z^{N-3}}+\cdots +\frac{1}{(N-1)!}\right]+\frac{1}{N!}\right)=\arctan\left(\frac{\Im(\xi)}{\Re(\xi)}\right),\end{eqnarray*}  
where  \begin{eqnarray*} \xi = \frac{1}{z}\left[\frac{1}{z^{N-1}}+\frac{1}{z^{N-2}}+\frac{1}{2!z^{N-3}}+\cdots +\frac{1}{(N-1)!}\right]+\frac{1}{N!}.\end{eqnarray*}Clearly $\Re(\xi)\geq \gamma_R$, where $\gamma_R \rightarrow 1/N!$ as $R \rightarrow \infty$, and $\Im(\xi)\leq M/R$ for some positive number $M$. Thus,\begin{eqnarray*} \arg\left(\frac{T_N(z)}{z^N}\right) =\arctan\left(\frac{\Im(\xi)}{\Re(\xi)}\right)\leq \arctan\left(\frac{M}{R\gamma_R}\right).\end{eqnarray*}Since $\arg\left(T_N(z)/z^N\right)=\arg\left(T_N(z)\right)-N\arg(z)-2q\pi$, for some integer $q$, we conclude from the above equation that \begin{eqnarray} \arg\left(T_N(z)\right)=N\arg(z)+2q\pi+\delta_R,\label{eq 5.9}\end{eqnarray}where $|\delta_R|\leq \arctan\left(M/R\gamma_R\right)$. Part (iii) now follows from $\left(\ref{eq 5.8}\right)$   and  $\left(\ref{eq 5.9}\right)$.
\end{proof}

\begin{lem} The function $e^z-T_N(z)$ has infinitely many zeros.  Furthermore, each of the nontrivial zeros is simple.\end{lem}

\begin{proof}  Assume on the contrary that $e^z-T_N(z)$ has a finite number of nontrivial zeros (or possibly none).  Let$z_1,\;z_2,\; \cdots, z_n$   be the nontrivial zeros and define $Q(z)=\prod_{k=1}^n(z-z_k)$  (or $Q(z)=1$ if there are no nontrivial zeros). By Weierstrass Factorization Theorem, we can express $e^z-T_N(z)$   as 
\begin{eqnarray}e^z-T_N(z)=P(z)e^{g(z)},\label{eq 5.10}\end{eqnarray}where $P(z)=z^{N+1}Q(z)/(N+1)!$. By comparing the growth rate of the two sides of $\left(\ref{eq 5.10}\right)$ (see \cite{HB}) we conclude that $g(z)=a z +b$ with $a\ne 0$. We now differentiate $\left(\ref{eq 5.10}\right)$ to get\begin{eqnarray}e^z-T_{N-1}(z)=\left(P'(z)+a P(z)\right)e^{az+b}.\label{eq 5.11}\end{eqnarray} Subtracting  $\left(\ref{eq 5.10}\right)$ from $\left(\ref{eq 5.11}\right)$ and noting that $T_N(z)-T_{N-1}(z)=z^N/N!$, we get\begin{eqnarray*}\frac{z^N}{N!}=\left((1-a)P(z)- P'(z)\right)e^{az+b}.\end{eqnarray*}  This last equation implies that $e^{az+b}$  is a rational function.  This contradiction shows that there are infinity many nontrivial roots of $\left(\ref{eq 5.1}\right)$.\\
\indent To prove the second statement, suppose to the contrary that $\omega$ is a root of multiplicity $m>1$. Then $ e^z-T_N(z)=(z-\omega)^m F(z)$, where $F(\omega)\ne 0$. As above, subtract the derivative of this last equation from itself to get\begin{eqnarray*}\frac{z^N}{N!}=(z-\omega)^{m-1}\left[(\omega-m-z)F(z)-(z-\omega) F'(z)\right].\end{eqnarray*}  However, the right hand side vanishes at $z=\omega$, while left hand side does not. Thus $m=1$  and the lemma follows.
\end{proof}

Next we specialize to the cases $N=2$ and $N=3$.  
\begin{lem}There is exactly one root $z_k=x_k+i y_k$ of $e^z-1-z=0$   having imaginary part inside the interval\begin{eqnarray} (2k+1/4)\pi<y_k<(2k+1/2)\pi\label{eq 5.12}\end{eqnarray} for each positive integer $k$  and no others besides their conjugates and $z=0$. \end{lem}
\begin{proof} Let $z=x+iy$ be a root with $y>0$.   Then equating real and imaginary parts of $e^z-1-z=0$, we get
\begin{eqnarray}
e^x\cos y -1 -x=0,\label{eq 5.13}\\e^x \sin y-y=0.\label{eq 5.14}
\end{eqnarray}
Since $e^x=y/\sin y>0$ from $\left(\ref{eq 5.14}\right)$, it follows that $2k\pi<y<(2k+1)\pi$ for some nonnegative integer $k$. But if $(2k+1/2)\pi<y<(2k+1)\pi$, then we must have $\cos y<0$, which forces $x<-1$  because of $\left(\ref{eq 5.13}\right)$.  It follows from $\left(\ref{eq 5.14}\right)$ that \begin{eqnarray*}y=e^x \sin y <e^x<\frac{1}{e}.\end{eqnarray*}   This contradicts the fact that $y>(2k+1/2)\pi\geq\pi/2$.  Hence,
\begin{eqnarray} 2k\pi<y_k<(2k+1/2)\pi.\label{eq 5.15}\end{eqnarray}	 

We now show that there is precisely one of $e^z-1-z=0$  satisfying  $\left(\ref{eq 5.12}\right)$.  First, we note that any solution of $\left(\ref{eq 5.13}\right)$ and $\left(\ref{eq 5.14}\right)$ must have $y$ as a root of \begin{eqnarray} f(y):=-1+y \cot y -\log\left(\frac{y}{\sin y}\right).\label{eq 5.16}\end{eqnarray}Then for $y$ satisfying $\left(\ref{eq 5.15}\right)$, we have \begin{eqnarray*} f'(y)=2 \cot y -y \csc^2 y-\frac{1}{y}=\frac{\sin(2y)-y}{\sin^2 y}-\frac{1}{y}<0.\label{eq 5.16}\end{eqnarray*} Hence $f(y)$ is strictly decreasing in interval $\left(\ref{eq 5.15}\right)$. Since $f\left((2k+1/4)\pi\right)>0$ and $f\left((2k+1/2)\pi\right)<0$, the lemma follows.
\end{proof}

\begin{lem}For each positive integer $k$, there is exactly one root $z_k=x_k+i y_k$ of $e^z-1-z-z^2/2=0$  such that  \begin{eqnarray} (2k+1/2)\pi<y_k<(2k+1)\pi.\label{eq 5.17}\end{eqnarray} Furthermore, there are no other roots besides their conjugates and $z=0$. \end{lem}

\begin{proof} We first prove existence.  As in the proof of Lemma  5.3, we let $z=x+iy$ be a root with $y>0$ and equate the real and imaginary parts of  $e^z-1-z-z^2/2=0$  to obtain 
\begin{eqnarray}e^x\cos y -1 -x-\frac{1}{2}\left(x^2-y^2\right)=0,\label{eq 5.18}\\e^x \sin y-y-x y=0.\label{eq 5.19}
\end{eqnarray}
By dividing the two equations above, we obtain \begin{eqnarray*}\cot y=\frac{1+x+\left(x^2-y^2\right)/2}{y(1+x)}.\end{eqnarray*}
The equation above is quadratic in $x$ and admits the solution set  \begin{eqnarray}x=-1+\frac{y}{\sin y}\left(\cos y\pm \sqrt{1-\frac{\sin ^2 y}{y^2}}\right).\label{eq 5.20}\end{eqnarray} We observe that the negative solution in $\left(\ref{eq 5.20}\right)$ is not allowed since this would imply from $\left(\ref{eq 5.19}\right)$ that\begin{eqnarray*}e^x=\frac{y(1+x)}{\sin y}=\frac{y^2}{\sin^2y}\left(\cos y-\sqrt{1-\frac{\sin ^2 y}{y^2}}\right)<\frac{y^2}{\sin^2y}\left(\cos y-\sqrt{1-\sin^2y}\right)=0,\end{eqnarray*}	 
which is a contradiction.  We can therefore rewrite $\left(\ref{eq 5.20}\right)$  as
\begin{eqnarray}1+x=\frac{y}{\sin y}\left(\cos y+\sqrt{1-\frac{\sin ^2 y}{y^2}}\right).\label{eq 5.21}\end{eqnarray}We also rewrite $\left(\ref{eq 5.19}\right)$ as \begin{eqnarray}e^x=\frac{y(1+x)}{\sin y}.\label{eq 5.22}\end{eqnarray}
We now take the natural log of $\left(\ref{eq 5.22}\right)$ and use $\left(\ref{eq 5.21}\right)$ to get \begin{eqnarray}x=\log \left[ \frac{y(1+x)}{\sin y} \right]= \log\left[ \frac{y^2} {\sin^2 y} \left(\cos y + \sqrt{ 1- \frac{ \sin^2 y}{y^2}}\right) \right]. \label{eq 5.23}\end{eqnarray}	

It follows from $\left(\ref{eq 5.21}\right)$ and $\left(\ref{eq 5.23}\right)$ that
\begin{eqnarray} -1+\frac{y}{\sin y}\left(\cos y + \sqrt{1-\frac{\sin ^2 y}{y^2}}\right)
= \log\left[ \frac{y^2} {\sin^2 y} \left(\cos y + \sqrt{1- \frac{ \sin^2 y}{y^2}}\right) \right].\label{eq 5.24}\end{eqnarray}	 
To show that $\left(\ref{eq 5.24}\right)$ admits a solution, we consider the function
	\begin{eqnarray} F(y)=-\sin y+y\left(\cos y + \sqrt{1-\frac{\sin ^2 y}{y^2}}\right)
- \sin y \log\left[ \frac{y^2} {\sin^2 y} \left(\cos y + \sqrt{ 1- \frac{ \sin^2 y}{y^2}}\right) \right].\label{5.25}\end{eqnarray} 	
The zeros of $F(y)$, excluding possibly those that are integer multiples of $\pi$, must be solutions of $\left(\ref{eq 5.24}\right)$.  We now apply the Intermediate Value Theorem to $F(y)$ to isolate these zeros.  To this end, we first note that for all positive integers $k$,
	  \begin{eqnarray} F((2k+1)\pi)=0& \mbox{and} & F'((2k+1)\pi)=1+\log\left[1/2(2k+1)^2\pi^2-1)\right]>0.\label{eq 5.26}\end{eqnarray}
Moreover, \begin{eqnarray} F((2k+1/2)\pi)=-1+\sqrt{(2k+1/2)^2\pi^2-1}-\log\left[(2k+1/2)^2\pi^2\sqrt{(2k+1/2)^2\pi^2-1}\right]>0.\label{eq 5.27}\end{eqnarray}

  Therefore, $F(y)$ possesses a root between $(2k+1/2)\pi$  and $(2k+1)\pi$  for all positive integers $k$  because of $\left(\ref{eq 5.26}\right)$ and $\left(\ref{eq 5.27}\right)$, and hence must be a root of $\left(\ref{eq 5.23}\right)$.  This completes the proof of the existence of a root with imaginary part in the desired interval.\\

We now prove uniqueness. To this end, set $I_k=((2k+1/2)\pi,(2k+1)\pi)$  and let $x=f(y)$ and $x=g(y)$  represent functions defined implicitly by $\left(\ref{eq 5.18}\right)$ and $\left(\ref{eq 5.19}\right)$, respectively.  We then differentiate implicitly to obtain
\begin{eqnarray}\frac{df}{dy}=\frac{e^{f(y)} \sin y - y}{e^{f(y)} \cos y -1 - f(y)}\label{eq 5.28},\\
\frac{dg}{dy}=-\frac{e^{g(y)} \cos y - 1-g(y)}{e^{g(y)} \sin y -y}.\label{eq 5.29}
\end{eqnarray}
Observe that $df/dy$  and $dg/dy$  are negative reciprocals of each other at every point of intersection between $f$ and $g$, i.e. whenever $f(y)=z=g(y)$, which  is nothing more than a restatement of the Cauchy-Riemann equations satisfied by $f$ and $g$.  Now, since $x>0$ and $\cos y<0$  for $y\in I_k$  , it follows from $\left(\ref{eq 5.28}\right)$  and $\left(\ref{eq 5.29}\right)$ that at each point of intersection we must have
	 \begin{eqnarray}\frac{df}{dy}=\frac{x y}{e^x \cos y - 1-x}<0,\label{eq 5.30}\\\frac{dg}{dy}=\frac{-e^x \cos y + 1+x}{x y}>0.\label{eq 5.31}\end{eqnarray}	

Next, we claim that $\left(\ref{eq 5.30}\right)$   and  $\left(\ref{eq 5.31}\right)$  restrict $f$ and $g$ to intersect at no more than one point inside $I_k$.  Assume on the contrary that they intersect at two distinct points with imaginary parts $y_1$  and $y_2$  with  $y_1<y_2$.  Then $\left(\ref{eq 5.30}\right)$ and $\left(\ref{eq 5.31}\right)$ imply that there exist $x$ and $y$ such that\begin{eqnarray*}y_1<x<y<y_2,\;\;g(x)>f(x),\;\;\mbox{and}\;\;g(y)<f(y).\end{eqnarray*}  
By the Intermediate Value Theorem, there exists a third intersection point at $y_3$, i.e. $f(y_3)=g(y_3)$, satisfying $y_1<y_3<y_2$.  Since $\left(\ref{eq 5.30}\right)$ and $\left(\ref{eq 5.31}\right)$ are also satisfied at $y_3$, we again have a fourth intersection point at $y_4$  satisfying $y_1<y_4<y_3$, and so forth.  This yields a bounded sequence of zeros   and therefore must contain an accumulation point.  It follows that the complex function $e^z-1-z-z^2/2$  is identically zero, which is a contradiction.  Hence, there is one and only one root with imaginary part inside $I_k$.
\end{proof}

\begin{rem} For $N = 2$, it follows from Lemma 5.2, $\left(\ref{eq 5.12}\right)$, and $\left(\ref{eq 5.14}\right)$ that the zeros of $e^z-1-z=0$  form a sequence of complex numbers with strictly increasing modulus. This assumption was made in deriving equation $\left(\ref{eq 4.8}\right)$. Similarly, for $N=3$, we use Lemma 5.2,  $\left(\ref{eq 5.17}\right)$ and $\left(\ref{eq 5.19}\right)$ to arrive at the same conclusion about the modulus of the roots. \end{rem}

\section{APPENDIX II}

\setcounter{equation}{0} 

Tables \ref{table:zerosN2} and \ref{table:zerosN3} list the first ten zeros of $e^z-T_{N-1}(z)=0$ for $N=2$ and $N=3$, respectively.  These values were computed using the software program Mathematica.

\begin{table}[p]
\begin{tabular}{|c|c|c|c|} \hline
$k$ & $z_k$ & $r_k$ & $\theta_k$ \\ \hline
1 & 2.088843016+7.461489286$i$ & 7.748360311 & 1.2978341024 \\ \hline
2 & 2.664068142+13.87905600$i$ & 14.13242564 & 1.3811541551 \\ \hline
3 & 3.026296956+20.22383500$i$ & 20.44900915 & 1.4222583654 \\ \hline
4 & 3.291678332+26.54323851$i$ & 26.74656346 & 1.4474143156 \\ \hline
5 & 3.501269010+32.85054823$i$ & 33.03660703 & 1.4646154233 \\ \hline
6 & 3.674505305+39.15107412$i$ & 39.32313052 & 1.4772159363 \\ \hline
7 & 3.822152869+45.44738491$i$ & 45.60782441 & 1.4868931567 \\ \hline
8 & 3.950805215+51.74088462$i$ & 51.89150222 & 1.4945866979 \\ \hline
9 & 4.064795694+58.03240938$i$ & 58.17459155 & 1.5008669923 \\ \hline
10 & 4.167125550+64.32248998$i$ & 64.45733203 & 1.5061018433 \\ \hline
\end{tabular} \\
\begin{tabular}{c}
\end{tabular}
\caption{First ten nonzero roots of $e^z-1-z=0$ in the upper-half complex plane.}
\label{table:zerosN2}
\end{table}

\begin{table}[p]
\begin{tabular}{|c|c|c|c|} \hline
$k$ & $z_k$ & $r_k$ & $\theta_k$ \\ \hline
1 & 3.838602048+8.366815507$i$ & 9.205349934 & 1.1406576364 \\ \hline
2 & 4.857263960+14.95891141$i$ & 15.72774757 & 1.2568294158 \\ \hline
3 & 5.520626554+21.39846201$i$ & 22.09912880 & 1.3183102795 \\ \hline
4 & 6.016178416+27.77895961$i$ & 28.42296607 & 1.3575169538 \\ \hline
5 & 6.412519686+34.12944500$i$ & 34.72663855 & 1.3850733959 \\ \hline
6 & 6.743013428+40.46233161$i$ & 41.02034263 & 1.4056646865 \\ \hline
7 & 7.026523305+46.78391852$i$ & 47.30863623 & 1.4217195916 \\ \hline
8 & 7.274789053+53.09777556$i$ & 53.59380865 & 1.4346366398 \\ \hline
9 & 7.495625078+59.40609018$i$ & 59.87710703 & 1.4452835555 \\ \hline
10 & 7.694499832+65.71028350$i$ & 66.15925246 & 1.4542298245 \\ \hline
\end{tabular} \\
\begin{tabular}{c}
\end{tabular}
\caption{First ten nonzero roots of $e^z-1-z-z^2/2=0$ in the upper-half complex plane.}
\label{table:zerosN3}
\end{table}

 \end{document}